\documentclass[preprint,12pt,authoryear]{elsarticle}

\usepackage{amssymb}
\usepackage{amsmath}
\usepackage{amsthm}

\usepackage{graphicx,subcaption}
\usepackage[dvipsnames]{xcolor}
\usepackage{makecell}
\usepackage[linesnumbered,ruled,vlined]{algorithm2e}
\usepackage{lscape}
\usepackage{rotating}
\usepackage{mathtools}
\usepackage{booktabs}
\usepackage{siunitx}
\usepackage{multirow}
\usepackage{lmodern}
\usepackage{array}

\usepackage{url}
\graphicspath{{JCP Submission/images/}}

\journal{Journal of Computational Physics}

\begin{document}

\begin{frontmatter}



\title{Neural Correction Operator: A Reliable and Fast Approach for Electrical Impedance Tomography} %% Article title


\author[1]{Amit Bhat\fnref{equal}} %% Author name
\author[2]{Ke Chen\fnref{equal}}
\author[1]{Chunmei Wang\fnref{equal}}
\fntext[equal]{These authors contributed equally to this work and are listed in alphabetical order.}

\affiliation[1]{organization={University of Florida},%Department and Organization
            city={Gainesville},
            postcode={32611}, 
            state={FL},
            country={}}

\affiliation[2]{organization={University of Delaware},%Department and Organization
            city={Newark},
            postcode={19716}, 
            state={DE},
            country={USA}}

%% Abstract
\begin{abstract}
Electrical Impedance Tomography (EIT) is a non-invasive medical imaging method that reconstructs electrical conductivity mediums from boundary voltage-current measurements, but its severe ill-posedness renders direct operator learning with neural networks unreliable. 
We propose the neural correction operator framework, which
learns the inverse map as a composition of two operators: a reconstruction operator using L-BFGS optimization with limited iterations to obtain an initial estimate from measurement data and a correction operator implemented with deep learning models to reconstruct the true media from this initial guess. We explore convolutional neural network architectures and conditional diffusion models as alternative choices for the correction operator.
We evaluate the neural correction operator by comparing with L-BFGS methods as well as neural operators and conditional diffusion models that directly learn the inverse map over several benchmark datasets.
Our numerical experiments demonstrate that our approach achieves significantly better reconstruction quality compared to both iterative methods and direct neural operator learning methods with the same architecture. 
The proposed framework also exhibits robustness to measurement noise while achieving substantial computational speedup compared to conventional methods.
The neural correction operator provides a general paradigm for approaching neural operator learning in severely ill-posed inverse problems. 
\end{abstract}

%% Keywords

\end{frontmatter}

\section{Introduction}
Electrical Impedance Tomography (EIT) is an inverse problem of finding the electrical conductivity distribution of an unknown medium via multiple voltage-current measurements on the domain boundary. 
Compared to typical methods like CT and MRI, EIT is a low cost, noninvasive, and radiation free method that has wide applications in medical imaging~\citep{meier2008assessment,conway1985applied,tarassenko1984imaging}, material engineering~\citep{duan20}, chemical engineering~\citep{waterfall1997visualizing,tapp2003chemical} and other fields.
The mathematical formulation of EIT is usually referred to as the Calder\'on problem~\citep{calderon2006inverse,uhlmann2009electrical} and its governing equation is the following elliptic partial differential equation (PDE).
\begin{equation}\label{eqn:elliptic}
\begin{aligned}
    -\text{div}\left(\sigma(x) \nabla u(x)\right) = 0\,, \quad x\in \Omega \,, \\
    u(x) = f(x)\,, \quad x \in \partial\Omega\,,
\end{aligned}
\end{equation}
where $\Omega\subset\mathbb{R}^d$ is a boundary Lipschitz domain, $u$ is the electrical potential distribution, $\sigma$ is the unknown conductivity distribution, and the Dirichlet boundary condition $f$ models the voltage applied on the boundary $\partial\Omega$.
The solution $u$ to \eqref{eqn:elliptic} is uniquely determined by $\sigma$ and $f$, as is the electrical current $g \coloneqq \sigma \frac{\partial u}{\partial n}|_{\partial\Omega}$, modeled by the Neumann derivatives measured on the boundary. 
A critical object in the Calder\'on problem is the Dirichlet-to-Neumann (D2N) map, defined as
\begin{equation}\label{eqn:D2Nmap}
    \Lambda_\sigma: H^{1/2}(\partial\Omega) \rightarrow H^{-1/2}(\partial\Omega)\,,\quad f \mapsto g \,.
\end{equation}
The Calder\'on problem consists of reconstructing $\sigma$ from the knowledge of the D2N map. The solution to the Calder\'on problem exists and can be uniquely determined under mild conditions~\citep{uhlmann2009electrical}. However, it is considered a severely ill-posed inverse problem due to poor stability~\citep{alessandrini1988stable,alessandrini1997examples}. 

In practice, measurement data are usually taken from finitely many receivers over a subset of the boundary, and are polluted with measurement noises. We denote the noisy dataset as $\mathcal{D}_N = \{(f_i,g_i)\ |\ g_i = \Lambda_\sigma f_i + \varepsilon_i\,,i=1,\ldots,N\} $, where the $f_i$'s and $g_i$'s denote the voltage and current measurements respectively. 
Reconstructing the medium $\sigma$ from partial data has been studied analytically in~\citet{bukhgeim2002recovering,sjostrand2004calderon,nachman2010reconstruction}; see~\citet{kenig2013recent} and the references therein for further details.
While EIT can be formulated in the Bayesian framework~\citep{dunlop2016bayesian}, numerical computation of the EIT problem from the measurement data $\mathcal{D}_N$ is usually formulated in the variational form~\citep{kohn1987relaxation,kohn1990numerical}. We consider the following variational form from~\citet{borcea2002electrical}: 
\begin{equation*}
    \min_\sigma \frac{1}{N} \sum_{i=1}^N\| \Lambda_\sigma f_i - g_i \|_{H^{-1/2}(\partial \Omega)}^2 \,.
\end{equation*}
In both the variational setting and Bayesian formulation, solving the EIT problem is computationally challenging. A large number of iterative solves of \eqref{eqn:elliptic} is usually involved, as the inverse map $\mathcal{D}_N\mapsto \sigma$ is typically numerically ill-conditioned. 
To mitigate this, various regularization terms are considered to encode prior knowledge of the target medium $\sigma$ under the optimization formulation, see~\citet{kaipio1999inverse,vauhkonen1998tikhonov}. In the Bayesian framework, such prior information is encoded via the prior distributions. 
However, these priors are handcrafted and therefore unable to accurately characterize the target medium $\sigma$, resulting in only marginal improvements and unsatisfactory reconstruction quality.

In the last decade, deep neural networks (DNNs) have achieved great success in computer vision, image processing and many machine learning tasks. 
More recently, their application to solving PDEs, in both forward and inverse settings, has become an emergent field of Scientific Machine Learning (SciML). 
Leveraging the universal approximation power of DNNs, various approaches have been developed to model unknown target functions, such as the deep image prior (DIP)~\citep{ulyanov2018deep}, physics-informed neural networks (PINN)~\citep{raissi2019physics,jagtap2022physics}, and many other works~\citep{pakravan2021solving,berg2017neural,lu2020extracting}.
Another strategy, known as \textit{operator learning}, focuses on using DNNs to directly learn the inverse map $\mathcal{D}_N \mapsto \sigma$ rather than modeling the medium function with a neural network surrogate. This strategy has been applied to several neural operators of various architectures~\citep{kovachki2023neural,molinaro2023neural,chen2024pseudo,wang2024latent,padmanabha2021solving,abhishek2024solving,cen2023electrical} as well as to conditional generative models like denoising diffusion probabilistic models (DDPMs)~\citep{ho2020denoising,song2021solving,chung2022diffusion,daras2024survey} and generative adversarial networks (GANs)~\citep{adler2018deep,patel2022solution}.

Despite the success of operator learning in solving forward PDEs, its application to inverse problems has mainly been focused on the linear case. Direct application of these universal neural operators to learn the inverse operator that maps the measurement data $\mathcal{D}_N$ to the target image $\sigma$ usually leads to inferior reconstructions~\citep{chen2024pseudo,chen2023let}. 
In this paper, we propose the neural correction operator framework, which approximates the target operator $\mathcal{D}_N \mapsto \sigma$ by expressing it as a composition of two operators: a reconstruction operator $\mathcal{R}$ that maps $\mathcal{D}_N$ to a rough reconstructed image $\hat{\sigma}$, and a correction operator $\mathcal{C}$ that maps $\hat{\sigma}$ to the true medium $\sigma$. We use L-BFGS with constant initialization and a limited number of iterations as the reconstruction operator $\mathcal{R}$ and neural operators as the correction operator $\mathcal{C}$. 
We compare our models against L-BFGS solvers with a large number of iterations and neural operators that directly learn the target operator $\mathcal{D}_N \mapsto \sigma $ over several benchmark datasets, including indicator functions with circular supports and Shepp-Logan phantoms. Numerical results demonstrate that the neural correction operator strategy significantly outperforms the baseline models, even when using the same neural operators with the same number of training data.
Our implementation of the methods and experiments described in this paper can be found at \url{https://github.com/amitbhat31/neural-correction-operator}.

\subsection{Related work and contributions}
In~\citet{guo2022transformer}, the authors modified the classical Direct Sampling Method (DSM) into a formulation that resembles the transformer model and applied it to solving the EIT problem. However, their method is limited to reconstructing only the support of $\sigma$ due to the original idea of DSM. 
Our method is able to reconstruct both the media support and values simultaneously, and can be easily adapted to a wide range of deep learning models.
Another related work is DeepEIT~\citep{liu2023deepeit}, where they follow the DIP idea by parametrizing the medium $\sigma$ with a neural network. The reconstructed image for a given measurement data can thus be obtained by training the neural network. This idea has the advantage of reconstructing images without any training data. However, their method can only deal with simple media distributions and is only capable of outperforming classical methods like TV-regularization and $\ell_2$-regularization. 
In our work, we have demonstrated the effectiveness of our method on several challenging datasets where $\ell_2$ regularization with the L-BFGS method yields poor reconstructions.
In the work of~\citet{abhishek2024solving}, the authors proposed using a DeepONet to learn the operators mapping Neumann-to-Dirichlet measurement data to the targeted conductivities. However, similar to standard neural operator learning approaches like ResNet, the reconstructions in~\citet{abhishek2024solving} suffer from blurry boundaries and large errors due to the ill-posedness of the EIT inverse problem.

To the best of our knowledge, the closest works to ours are the following. 
In~\citet{wei2019dominant}, the authors use the bases-expansion subspace optimization method (BE-SOM) to obtain multiple polarization tensors, followed by a convolutional neural network (CNN) to learn the relationship between input channels of polarization tensors and output channels of reconstructed media. The final reconstruction is the average of the CNN output channels. 
However, as a consequence of BE-SOM, the CNN requires $n$ input channels, where $n$ is the total number of sources in the measurement data. Such a requirement leads to a large and computationally expensive neural network. In comparison, our method only requires an initial guess generated from the L-BFGS algorithm. 
A similar idea was considered in~\citet{chen2020electrical}, where the author aims to reconstruct the support of the conductivities. Specifically, a conditional GAN was used to learn the relationship between the initial guess of the medium and the target medium. Here, the initial guesses are generated via the Landweber and Newton-Raphson algorithms. However, the Newton-Raphson algorithm involves constructing the Hessian matrix, which is not practical for PDE inverse problems.
In comparison, our method adopts the L-BFGS algorithm to generate the initial guess, which does not involve calculating the Hessian and can still obtain superlinear convergence. Furthermore, we have applied the ResNet and conditional DDPM formulations to improve the initial guesses. As a consequence, our method is able to generate reconstructed images with significantly sharper boundaries, all while avoiding expensive computing to generate the initial guesses. 

Another line of related work uses U-Net as a post-processing step to improve initial guesses generated using direct reconstruction methods, including the Direct Sampling Methods~\citep{guo2021construct}, D-bar methods~\citep{hamilton2018deep}, and the Calder\'on method~\citep{cen2023electrical}. However, these direct methods are usually much more expensive than iterative methods~\citep{wei2019dominant,chen2020electrical}, thus leading to a long inference time. 

Our work can also be viewed as a deep learning--based postprocessing step for traditional L-BFGS methods. A related line of ``warm-initiated'' methods in~\citet{zhou2023neural} and~\citet{guo2025warm} instead uses deep learning models as a pre-processing stage for gradient-based solvers, with applications to inverse scattering and fluorescence molecular tomography, respectively. The key difference lies in the role of the neural operator: in our approach, it acts as an image-to-image correction operator, while in warm-initiated methods it maps measurement data to target images to provide a good initial guess for subsequent gradient-based iterations. This difference in roles also leads to different architectural designs. In particular, our method can use a standard U-Net because both input and output lie in the same parameter space, whereas warm-initiated methods require carefully designed encoders to transform measurement space into parameter space. In terms of computational cost, one disadvantage of our method compared with warm-initiated methods is that it requires additional computation cost in both offline data generation and online inference.

\subsection{Organization of the paper}
The rest of the paper is structured as follows. We provide essential background on L-BFGS and DDPM formulations in Section~\ref{sec:background}. Section~\ref{sec:method} introduces our proposed neural correction operator framework. Section~\ref{sec:experiments} discusses the experimental setup and presents our numerical results. Finally, Section~\ref{sec:conclusion} summarizes the conclusions and discusses the challenges and future directions arising from this work.

\section{Background} \label{sec:background}
\subsection{L-BFGS method} 
The limited-memory Broyden–Fletcher–Goldfarb–Shanno (L-BFGS) method is a quasi-Newton optimization algorithm designed for large-scale unconstrained problems $\min_{x\in \mathbb{R}^d} \mathcal{L}(x)$. Like the classical BFGS method, L-BFGS approximates the inverse Hessian matrix of $\mathcal{L}$ using only gradient evaluations. However, instead of storing a full $d\times d$ matrix, L-BFGS maintains a limited history of the most recent $m$ pairs of iterate and gradient differences, denoted by
\[
s_k = x_{k+1} - x_k, \quad y_k = \nabla \mathcal{L}(x_{k+1}) - \nabla \mathcal{L}(x_k),
\]
which are used to implicitly construct a low-rank approximation of the inverse Hessian and to compute the descent direction.

This limited-memory approach enables L-BFGS to scale efficiently in high-dimensional settings. The method is typically coupled with a line search satisfying the Wolfe conditions to ensure global convergence. We refer the reader to~\citet{zhu1997algorithm} for more details.

In this work, we use L-BFGS not only as an optimization tool, but also to extract initial guesses of the reconstructed media, which is shown to be beneficial for downstream inference and learning tasks.

\subsection{Denoising Diffusion Probabilistic Models}
Denoising diffusion probabilistic models (DDPMs) are a category of score-based generative models that learn a target data distribution from a dataset of its samples through a forward and a reverse denoising process~\citep{ho2020denoising}.
In the forward process, Gaussian noises at various scales are added to the samples from the target distribution until they approach an isotropic Gaussian distribution. 
The reverse process learns how to successively denoise a sample from the standard Gaussian distribution back to a sample from the target distribution. 
Below, we provide an overview of both unconditional and conditional DDPM.

\subsubsection{Unconditional DDPM}
For convenience, we adopt the stochastic differential equation (SDE) formulation of DDPM developed in~\citet{song2021sbgm}.
Let $p(x)$ be the target data distribution from which the dataset is constructed. For time $t \in [0, T]$, a general framework for the forward process of score-based generative models can be expressed by the solution to the following SDE:
\begin{equation}\label{eqn:uncforwardSDE}
    dx_t = f(x_t, t)dt + g(t)dw_t, \quad t \in [0, T],
\end{equation}
where $f(x_t, t): {\mathbb{R}^n} \to \mathbb{R}^n$ and $g(t): \mathbb{R} \to \mathbb{R}$ are functions called the drift and diffusion coefficients of $x_t$, respectively, and $w_t$ is a standard Brownian motion. 
We denote the marginal probability distribution of $x_t$ as $p_t(x_t)$ and the transition distribution from $x_s$ to $x_t$ as $p_{st}(x_t | x_s)$ for $0 \leq s < t \leq T$. 
Starting from samples $x_0 \sim p_0(x_0) \equiv p(x)$, noise is gradually added via \eqref{eqn:uncforwardSDE} to obtain samples $x_T \sim p_T(x_T)$, where $p_T(x_T)$ follows the standard Gaussian distribution. 

The reverse process aims to start from Gaussian samples $x_T \sim p_T(x_T)$ and gradually denoise them to recover target samples $x_0 \sim p(x)$. This process is described by the reverse-time SDE~\citep{anderson1982revdiff}:
\begin{equation}\label{eqn:uncreverseSDE}
    dx_t = [f(x_t, t) - g(t)^2\nabla_x\log p_t(x_t)]dt + g(t)d\bar{w}_t    
\end{equation}
where $\bar{w}_t$ is a backward Brownian motion and $\nabla_x \log p_t(x_t)$ is called the score function. 
Once the score function is known, solving \eqref{eqn:uncreverseSDE} allows us to generate target samples $x_0$ from Gaussian samples $x_T$.
In practice, $\nabla_x \log p_t(x_t)$ is approximated by a neural network $s_\theta(x_t, t)$, where $\theta$ denotes learnable parameters. 
The score function $s_\theta(x_t, t)$ can be learned through the following training loss~\citep{hyvarinen2005ISM,vincent2011connection},
\begin{equation}\label{eqn:uncELBOsde}
    \mathcal{L}(\theta) = \mathbb{E}_{x_0, t, x_t | x_0} [\| s_\theta(x_t, t) - \nabla_{x_t}\log p_{0t}(x_t | x_0)\|_2^2].
\end{equation}
Generating samples using \eqref{eqn:uncreverseSDE} via a pre-trained score function in \eqref{eqn:uncELBOsde} resembles sampling Langevin dynamics, from which convergence of the reverse process to $p(x)$ is guaranteed~\citep{lee2023conv}.

The DDPM formulation of score-based generative models discretizes the above process for $t = [1, \ldots, T]$ by way of a variance schedule $\{\beta_t\}_{t=1}^T$ such that $0 < \beta_1 < \beta_2 < \cdots \beta_T < 1$. The variance schedule describes how noise is added at each step. The discrete-time Markov chain is described by
\begin{equation}\label{eqn:uncDDPMtrans}
    x_t = \sqrt{\bar{\alpha}_t}x_0 + \sqrt{1 - \bar{\alpha}_t}\epsilon_t, \quad t = 1, \ldots, T, 
\end{equation}
where $\alpha_t = 1 - \beta_t$, $\bar{\alpha}_t = \prod_{i=1}^t \alpha_i$, and $\epsilon_t \sim \mathcal{N}(0, \mathbf{I})$~\citep{ho2020denoising}. As $T \to \infty$, \eqref{eqn:uncDDPMtrans} converges to \eqref{eqn:uncforwardSDE} with $f(x_t, t) = -\frac{1}{2}\beta_t x_t$ and $g(t) = \sqrt{\beta_t}$~\citep{song2021sbgm}. 
We can simplify the score-matching objective described in \eqref{eqn:uncELBOsde}, which optimizes $s_\theta(x_t, t)$ into the following denoising evidence lower bound (ELBO) loss, which optimizes the denoiser function $\epsilon_\theta(x_t, t)$ as described in~\citet{ho2020denoising}:
\begin{equation}\label{eqn:uncELBOddpm}
    \mathcal{L}_{\mathrm{DDPM}}(\theta) = {\mathbb{E}}_{x_0,t,\epsilon_t}\left[\frac{\beta_t}{2\alpha_t(1 - \bar{\alpha}_t)}\| \epsilon_t - \epsilon_\theta(x_t, t)\|_2^2\right].
\end{equation}
Consequently, instead of learning the score function approximation $s_\theta(x_t, t)$, a neural network is used to learn the denoiser $\epsilon_\theta(x_t, t)$.
Then, sampling via the reverse process in the discrete-time setting and in terms of $\epsilon_\theta(x_t, t)$ is:
\begin{equation}\label{eqn:uncDDPMrevtrans}
    x_{t-1} = \frac{1}{\sqrt{\alpha_t}}\left(x_t - \frac{\beta_t}{\sqrt{1 - \bar{\alpha}_t}}\epsilon_\theta(x_t, t)\right) + \sqrt{\beta_t}z_t, \quad t= T, \ldots, 1.
\end{equation}
where $z_t \sim \mathcal{N}(0, \mathbf{I})$.

\subsubsection{Conditional DDPM}\label{sec:condDDPM}

The unconditional DDPM framework aims to generate samples from the distribution $p_0(x_0)$ of the medium without additional knowledge. 
In contrast, solving the inverse problem~\eqref{eqn:elliptic} requires learning a medium $x_0$ given its corresponding measurement $y$. Thus, the goal is to learn the measurement-to-medium operator $y\mapsto x_0$ rather than find an arbitrary sample from the distribution of the true media.
To this end, conditional DDPM can be used to conduct operator learning, i.e., generating samples from the posterior distribution $p(x_0 \ |\ y)$.

Applying conditional DDPM to learning PDE operators has been applied in several works, including PDE downscaling~\citep{lu2024pgdm} and PDE-based data assimilation~\citep{shysheya2024conditional}. For computational simplicity, we adopt a data-driven approach similar to that of~\citet{lu2024pgdm} and simply treat the measurement $y$ as an additional input in our noise approximator $\epsilon_\theta(x_t, y, t)$.

As a result of this formulation, the forward process remains the same as described in \eqref{eqn:uncDDPMtrans}, and the ELBO training loss \eqref{eqn:uncELBOddpm} becomes
\begin{equation}\label{eqn:condELBOddpm}
    \mathcal{L}_{\mathrm{cond}}(\theta) = {\mathbb{E}}_{x_0,t,\epsilon_t}\left[\frac{\beta_t}{2\alpha_t(1 - \bar{\alpha}_t)}\| \epsilon_t - \epsilon_\theta(x_t, y, t)\|_2^2\right],
\end{equation}
for $\epsilon_t \sim \mathcal{N}(0, \mathbf{I})$ and $t \sim \mathcal{U}(\{1, \ldots, T\})$.
Similarly, the reverse process is:
\begin{equation}\label{eqn:confDDPMrevtrans}
    x_{t-1} = \frac{1}{\sqrt{\alpha_t}}\left(x_t - \frac{\beta_t}{\sqrt{1 - \bar{\alpha}_t}}\epsilon_\theta(x_t,y, t)\right) + \sqrt{\beta_t}z_t, \quad t=T, \ldots, 1.
\end{equation}
where $z_t \sim \mathcal{N}(0, \mathbf{I})$.
Details on training and sampling for conditional DDPMs in the context of the EIT problem are provided in Section~\ref{sec:cDDPM}.

\section{Neural Correction Operator}\label{sec:method}

We first discretize the computational domain $\Omega \subset \mathbb{R}^2$ via a finite element mesh with $n_e$ elements and $n_b$ boundary nodes. 
The conductivity medium is then discretized as $\sigma \in \mathbb{R}^{n_e}$ and its corresponding D2N measurement can be represented as a matrix $M \in \mathbb{R}^{n_b \times n_b}$. 
Our objective is to learn the inverse operator
\begin{align*}
\mathcal{F}:\mathbb{R}^{n_b \times n_b} \rightarrow \mathbb{R}^{n_e}, \quad \mathcal{F}\left(M\right) = \sigma.
\end{align*}
We propose the neural correction operator, which decomposes $\mathcal{F}$ as the composition of two operators: 
\begin{equation}\label{eqn:method}
\begin{aligned}
    \mathcal{F} = \mathcal{C} \circ \mathcal{R}_K.
\end{aligned}  
\end{equation}
Here, the reconstruction operator $\mathcal{R}_K: \mathbb{R}^{n_b \times n_b} \to \mathbb{R}^{n_e}$ maps $M$ to a low-fidelity reconstruction $\hat{\sigma}$ via the L-BFGS method with fixed initialization and a fixed number $K$ of iterations, while the correction operator $\mathcal{C}: \mathbb{R}^{n_e} \to \mathbb{R}^{n_e}$ maps $\hat{\sigma}$ to the true medium $\sigma$. The initialization is chosen as a constant function with values equal to the known background of the conductivity medium. The value of $K$ is usually small and depends on the medium distribution and noise level. More details on the choice of $K$ will be discussed in Section~\ref{sec:setup}.

As the reconstruction operator $\mathcal{R}_K$ is determined, learning the target operator $\mathcal{F}$ is now reduced to learning the correction operator $\mathcal{C}$.
When noises are small, we propose to use a deep learning model, e.g., ResNet, to learn the correction operator. 
In the high-noise regime, the posterior distribution $p(x_0\ |\ y)$ can be multi-modal. Therefore, the inverse operator $\mathcal{F}$ as well as the correction operator $\mathcal{C}$ may not be well-defined and can be multi-valued. To this end, we propose to use a conditional diffusion model to approximate the multi-valued operator $\mathcal{C}$.

\subsection{Reconstruction Operator}

\begin{figure}[t]
\centering
\begin{subfigure}[t]{0.49\linewidth}  
    \centering 
    \includegraphics[width=\linewidth]{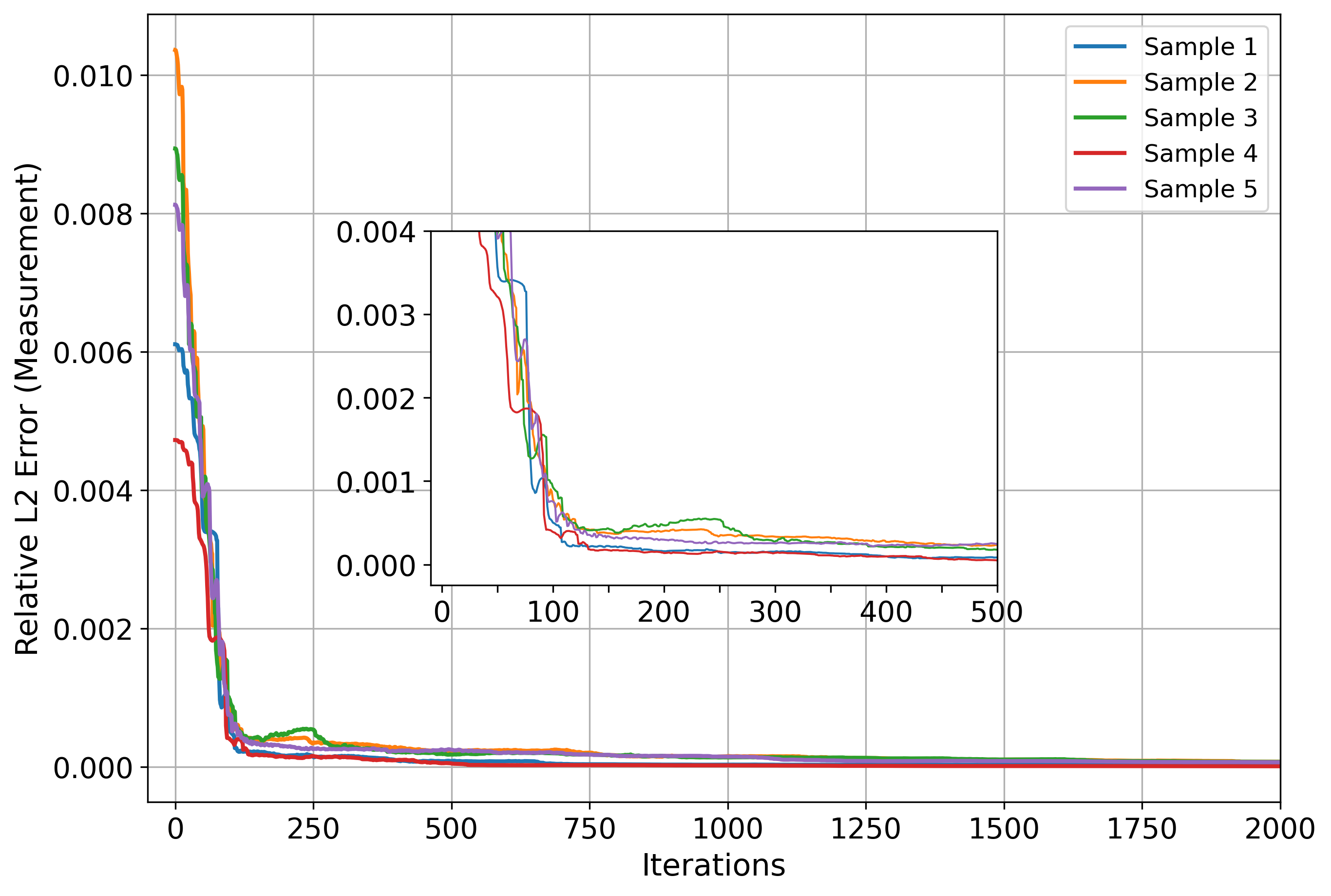}
    \caption[]%
    {No noise case.}  
    \label{fig:no_noise_A}
\end{subfigure}
\hfill
\begin{subfigure}[t]{0.49\linewidth}
    \centering
    \includegraphics[width=\linewidth]{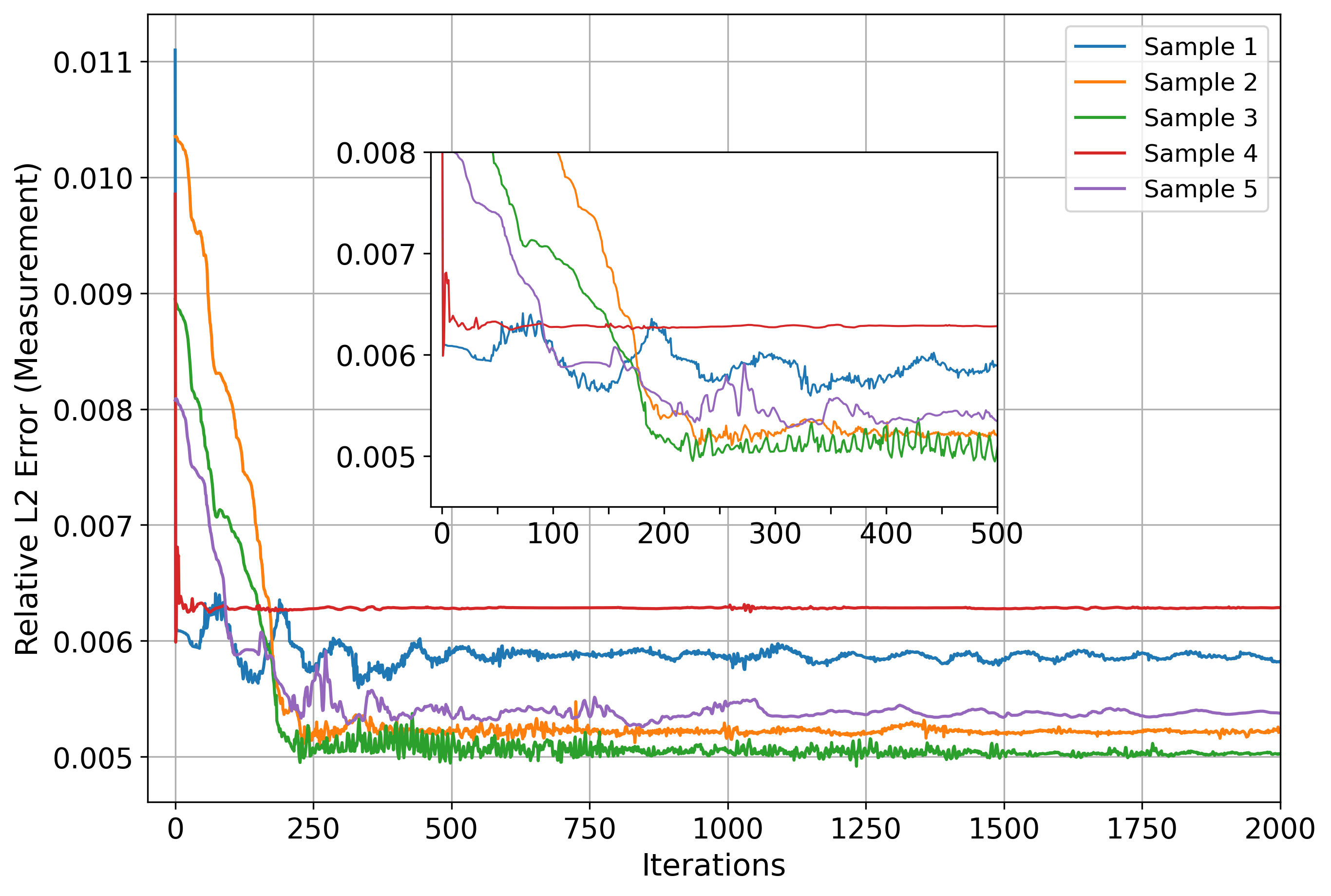}
    \caption[]%
    {1\% noise case.}
    \label{fig:1pct_noise_B}
\end{subfigure}
\caption{Convergence of L-BFGS for Shepp-Logan media.} 
\label{fig:conv_lbfgs}
\end{figure}

Conventional nonlinear optimization methods, such as L-BFGS, perform poorly in solving PDE inverse problems due to the ill-posedness and non-convexity.
In particular, reconstructed images usually suffer from blurry edges and missing finer details.
For a given D2N measurement $M \in \mathbb{R}^{n_b \times n_b}$, we consider the following optimization problem
\begin{equation}\label{eqn:rec}
\underset{\sigma}{\mathrm{argmin}} ||M - M_\sigma||_2^2 \,,
\end{equation}
where $M_\sigma$ denotes the D2N measurement generated by $\sigma$. 
We now specify the definition of $\mathcal{R}_K: M \mapsto \hat{\sigma}$.
We consider the L-BFGS solver to \eqref{eqn:rec}, initiated with a constant iterate $\sigma_{init} \in \mathbb{R}^{n_e}$ with the known background value, and we define $\hat{\sigma}$ as the $K$-iterate of the L-BFGS solver. 
Due to the non-convexity of \eqref{eqn:rec}, $\hat{\sigma}$ generated from a large number of iterations is not necessarily close to a global minimum. 
The ill-posedness of EIT also prevents more iterations of L-BFGS from further improving the reconstruction quality. We highlight this behavior in Figure~\ref{fig:conv_lbfgs}, which display the convergence of L-BFGS on five ground truth Shepp-Logan media in noiseless and 1\% noisy conditions. From Figure~\ref{fig:no_noise_A}, we observe that the relative $\ell_2$ error of the measurement of the L-BFGS-reconstructed media saturates after $K = 150$ L-BFGS iterations when there is no noise in the measurement data. In the $1\%$ noise case, we found that $K = 250$ iterations were necessary for the error to saturate, which can be seen in Figure~\ref{fig:1pct_noise_B}. With negligible improvement in the reconstructed measurement after these values of $K$, we choose these values in order to achieve the best tradeoff between initial guess quality and computational cost. 
By keeping the number of iterations relatively low, we generate a low-fidelity reconstruction $\hat{\sigma}$ with fast evaluation of $\mathcal{R}_K$.

\subsection{Approximation of the Correction Operator}
To learn the correction operator $\mathcal{C}$, we construct a new training dataset $\hat{\mathcal{D}}_N \coloneqq \{(\hat{\sigma}^{(i)},\sigma^{(i)}),\, i=1,\ldots,N\}$, where each $\hat{\sigma}^{(i)}$ is computed via $\hat{\sigma}^{(i)} \coloneqq \mathcal{R}_K(M^{(i)})$. This process requires offline computation, comprising of a small number of L-BFGS solves applied to the original training dataset $\mathcal{D}_N = \{(M^{(i)},\sigma^{(i)}), \allowbreak \, i=1,\ldots,N\}$.
Then, a deep learning model is used to learn the correction operator $\mathcal{C}$ from the new dataset $\hat{\mathcal{D}}_N$. We explore ResNet and conditional DDPMs as alternatives for the correction operator, and discuss our formulations of these models in the following subsections.
For compatibility with these models, we interpolate $\sigma$ and $\hat{\sigma}$ into square images, i.e.,  $\sigma, \hat{\sigma} \in \mathbb{R}^{n_i \times n_i}$ where $n_i$ is the chosen image size.

\subsubsection{ResNet}\label{sec:ResNet}
ResNet is a discriminative DNN architecture with skip connections that is introduced in~\citet{he2016deep} to address the vanishing gradient issues in training deep neural networks. 
The standard ResNet architecture consists of convolutional layers that downsample the inputs at each step to reduce spatial resolution and aggregate features for the purpose of classification.  

\begin{figure}[t]
    \centering
    \includegraphics[width=0.95\linewidth]{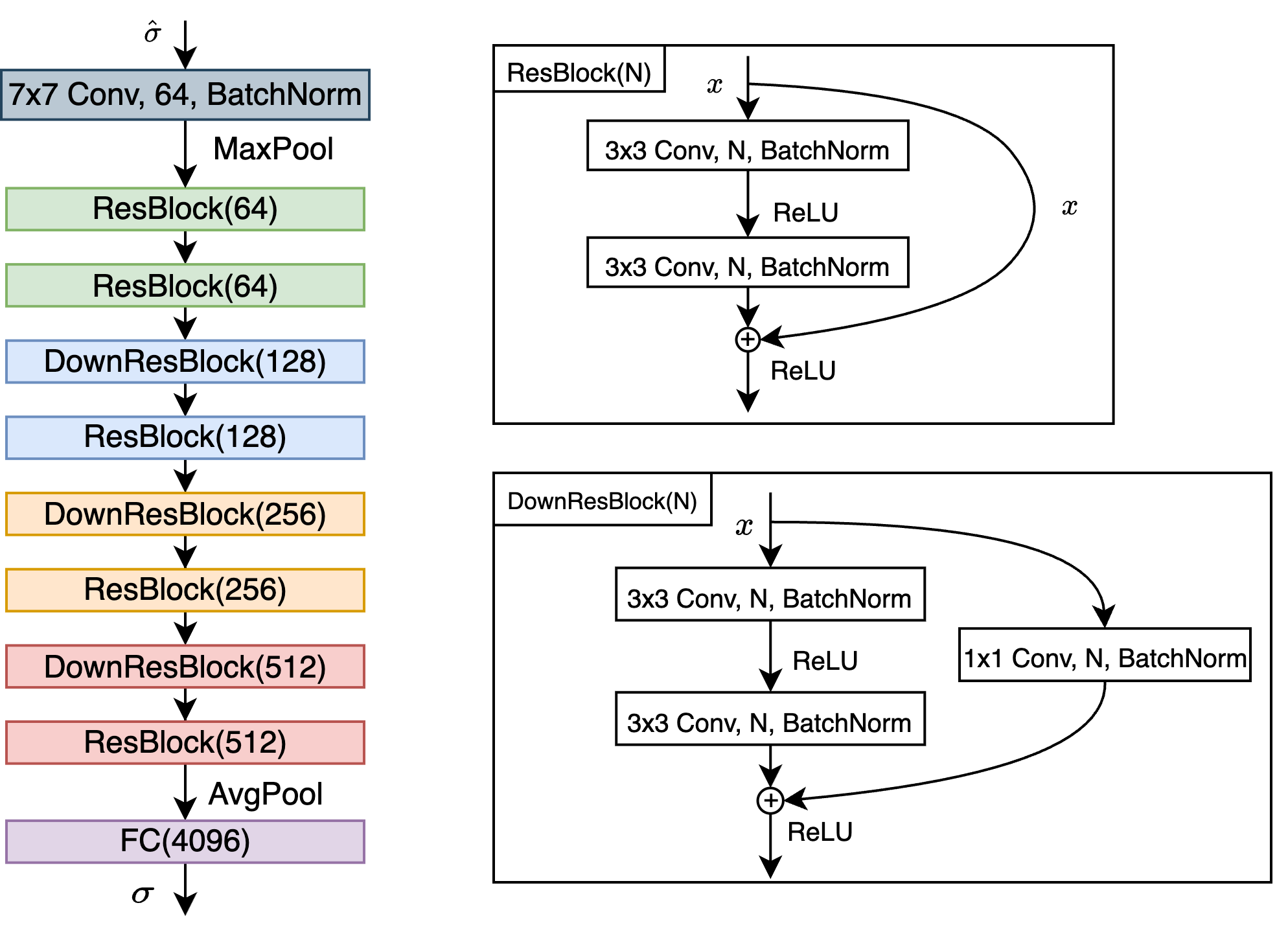}
    \caption{\textbf{Left:} ResNet architecture used to learn the neural operator $\mathcal{C}_R$. We use 8 residual blocks to learn the overall features of the image, and use a fully-connected layer at the end to upsample back to the input dimensions. \textbf{Right:} Composition of a ResBlock and a DownResBlock. $N$ denotes the number of channels in the input.}
    \label{fig:resnet-arch}
\end{figure}

We approximate the target operator $\mathcal{C}$ with a ResNet model $\mathcal{C}_R$ with its architecture shown in Figure~\ref{fig:resnet-arch}. 
As we are dealing with an image reconstruction task, we replace the fully-connected final layer in typical ResNet implementations with one that acts as an upsampling layer back to the original resolution of the image. We employ this fully-connected upsampling layer to reduce the number of parameters of our model as opposed to using a typical convolutional upsampling approach. 
We train $\mathcal{C}_R$ with mean squared error (MSE) loss $\mathcal{L}_\text{MSE}$:
\begin{equation}
    \mathcal{L}_\text{MSE} = \frac{1}{M}\sum_{i=1}^M\|\sigma - \mathcal{C}_R (\hat{\sigma})\|_2^2.
\end{equation}
Further training details such as the choice of optimizer and learning rates are discussed in Section~\ref{sec:traindets}.

\subsubsection{Conditional DDPM}\label{sec:cDDPM}
We utilize a conditional DDPM model as another method to learn $\mathcal{C}$ to generate samples of $\sigma$ conditioned on the initial guess $\hat{\sigma}$. 
Given $\sigma_0 = \sigma$, the forward process is given by the discrete-time Markov chain described in~\eqref{eqn:uncDDPMtrans}. 
Similarly, the learned reverse process starts from $\sigma_T \sim \mathcal{N}(0, \mathbf{I})$ and is given by~\eqref{eqn:confDDPMrevtrans}. 
We train the denoiser neural network $\epsilon_\theta(\sigma_t, \hat{\sigma}, t)$ on the simplified denoising ELBO loss as discussed in~\citet{ho2020denoising}:
\begin{equation*}
    \mathcal{L}_{\mathrm{simple}}(\theta) = {\mathbb{E}}_{\sigma_0,t,\epsilon_t}\left[\| \epsilon_t - \epsilon_\theta(\sigma_t, \hat{\sigma}, t)\|_2^2\right].
\end{equation*}
From these formulations, we describe the pretraining and sampling algorithms in Algorithms~\ref{alg:train} and~\ref{alg:samp} respectively.

\begin{algorithm}[t]
\caption{Training phase}\label{alg:train}
\SetAlgoLined
\LinesNumbered
\Repeat{$\mathrm{converged}$}{
    $\sigma = \sigma_0 \sim p(\sigma), \, \sigma \in \mathbb{R}^{n_i \times n_i}$\;
    $\hat{\sigma} = \mathcal{R}_K(M), \, \hat{\sigma} \in \mathbb{R}^{n_i \times n_i}$\;
    $t \sim \mathcal{U}(\{1, \ldots, T\})$\;
    $\epsilon_t \sim \mathcal{N}(0, \mathbf{I}), \, \epsilon_t \in \mathbb{R}^{n_i \times n_i}$\;
    Perform gradient descent step on:
    $\nabla_\theta ||\epsilon_t - \epsilon_\theta(\sqrt{\bar{\alpha}_t}\sigma_0 + \sqrt{1 - \bar{\alpha}_t}\epsilon_t, \hat{\sigma}, t)\|_2^2$
    }
\end{algorithm}

\begin{algorithm}[t]
\caption{Sampling phase}\label{alg:samp}
\SetAlgoLined
\LinesNumbered
$\sigma_T \sim \mathcal{N}(0, \mathbf{I}), \sigma_T \in \mathbb{R}^{n_i \times n_i}$\;
\For{$t = T, \ldots, 1$} {
    $z \sim \mathcal{N}(0, \mathbf{I})$ if $t > 1$ else $z = 0$\;
    $\sigma_{t-1} = \frac{1}{\sqrt{\alpha_t}}\left( \sigma_t - \frac{\beta_t}{\sqrt{1 - \bar{\alpha}_t}}\epsilon_\theta(\sigma_t, \hat{\sigma}, t)\right) + \sqrt{\beta_t}z$\;
}
\Return{$\sigma_0$}\;
\end{algorithm}

\begin{figure}[t]
    \centering
    \includegraphics[width=0.9\linewidth]{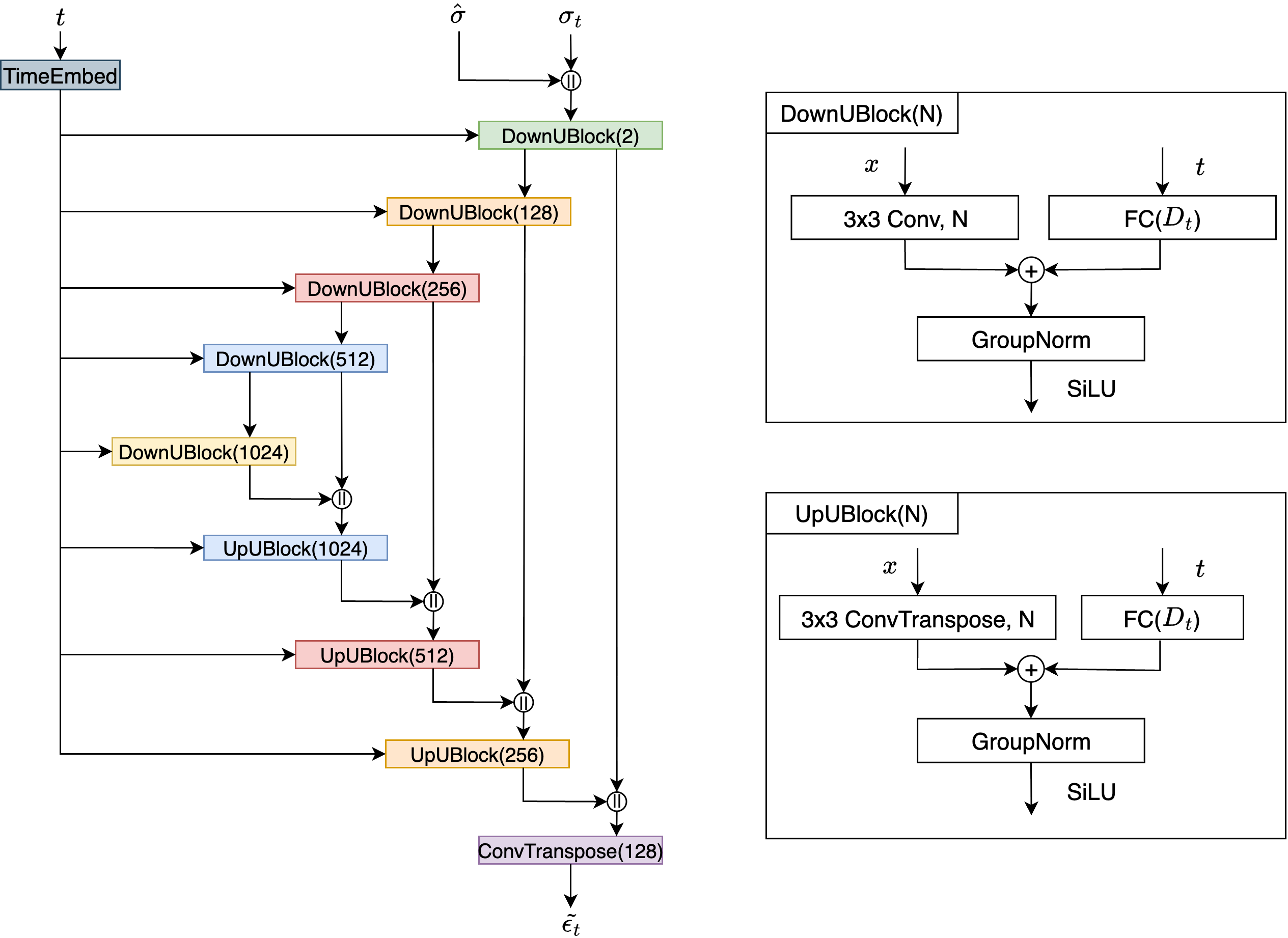}
    \caption{\textbf{Left:} UNet architecture used to learn the neural operator $\mathcal{C}_R$. We use 4 downsampling blocks to learn the overall features of the image, and 4 upsampling with residual connections to return to the input dimensions. Here ``$||$'' denotes concatenation along the channel dimension and ``$+$'' denotes addition. 
    \textbf{Right:} Composition of a DownUBlock and an UpUBlock. Here $N$ denotes the number of channels in the input and $D_t$ is the embedding dimension.}
    \label{fig:unet-arch}
\end{figure}

In our DDPM implementation, a UNet without attention is used to learn the denoiser $\epsilon_\theta(\sigma_t, \hat{\sigma}, t)$ as proposed by~\citep{ho2020denoising, song2021solving}. Figure~\ref{fig:unet-arch} illustrates our specific architecture. The UNet takes as inputs the noised image $\sigma_t$, the initial guess $\hat{\sigma}$, and the timestep $t$ and outputs a noise prediction $\tilde{\epsilon_t}$. We use sinusoidal embeddings followed by a fully-connected layer in the initial TimeEmbed block to yield a time embedding vector of size $D_t = 256$. 
Further training details such as the choice of optimizer and learning rates are discussed in Section~\ref{sec:traindets}.
Once we obtain a converged $\epsilon_\theta$, we can tractably sample and obtain an approximation to $\sigma$ via the sampling procedure in Algorithm~\ref{alg:samp}.

\section{Numerical Experiments}\label{sec:experiments}

\subsection{Experimental Setup}\label{sec:setup}
For all experiments, we assume the computational domain $\Omega$ is the unit disk $D_1$. 
We discretize $\Omega$ using $n_e = 2774$ triangular elements with $n_b = 128$ boundary nodal points. 
We compare our neural correction operator methods with several baseline models such as L-BFGS and neural operator methods for the EIT inverse problem over two benchmark datasets:
\begin{itemize}
    \item \textbf{Four Circles Distribution:} 
    The medium consists of a unit background value and several indicator functions with circular supports within $\Omega$, defined as the following. 
    \begin{equation*}
        \sigma_S(x) = 1 + \sum_{i \in S} w_i \cdot \mathbf{1}_{\{\|x - c_i\| \leq r_i\}}, \quad x \in \Omega, 
    \end{equation*}
    where the constant $1$ term denotes the background value of the media, and $c_i\in\mathbb{R}^2$ and $r_i\in\mathbb{R}_+$ denotes the random center and radius of the circular supports respectively.
    In particular, we select a subset $S \subseteq [4]$ uniformly at random to determine whether a particular random circular indicator will be included. 
    For $i \in [4]$, the centers $c_i \in \Omega$ and radii $r_i$ are independently sampled as the following:
    \begin{align*}
        c_1 &\sim \mathcal{U}([0.1, 0.4]^2), \\
        c_2 &\sim \mathcal{U}([-0.4, -0.1] \times [0.1, 0.4]), \\
        c_3 &\sim \mathcal{U}([-0.4, -0.1]^2), \\
        c_4 &\sim \mathcal{U}([0.1, 0.4] \times [-0.4, -0.1]), \\
        r_i &\sim \mathcal{U}(0.1, 0.4), \quad i \in [4],
    \end{align*}
    and height values $w_i$ of each circle are defined as $w_i = 2i$ for $i \in [4]$. Reconstructing this media distribution can be challenging as the contrast of the media can be as high as $9:1$.
    
    \item \textbf{Shepp-Logan Phantom Distribution:} Shepp-Logan phantoms are a commonly used test distribution in medical imaging, serving as a model of a human head~\citep{shepplogan1989}. They are defined as indicator functions supported on ellipses. For our experiments, we randomly vary the axis lengths, positions, and rotation angles to generate a diverse dataset~\citep{ruthotto2018slIP}. 
\end{itemize}
Unless otherwise established, we employ consistent data generation protocols for the neural correction operator methods. For each dataset, we create $N = 5000$ data pairs $\{\sigma_i, \hat{\sigma}_i\}_{i=1}^{N}$ via the following procedure: 
\begin{enumerate}
    \item Generate the true image $\sigma_i$ over the unit disk $D_1$.
    \item Compute the corresponding D2N measurement $M_i$ by solving \eqref{eqn:elliptic} with the finite element method (FEM). 
    \item Compute the low-fidelity reconstruction $\hat{\sigma}_i$ by solving \eqref{eqn:rec} with L-BFGS run for $K$ iterations with initial iterate $\sigma_{init} = \mathbf{1} \in \mathbb{R}^{n_e}$, where $K=350$ for the Four Circles distribution and $K=150$ for the Shepp-Logan distribution. 
\end{enumerate}
To facilitate the use of deep learning models, each medium $\sigma_i$ is resampled onto a uniform grid over $[-1, 1]^2$ using linear interpolation and constant padding with value 1; we denote these converted media as $\sigma_i^s$.
The resulting dataset $\{\sigma_i^s, \hat{\sigma}_i^s\}_{i=1}^{N}$ is randomly partitioned into 4,000 training pairs, 100 validation pairs, and 900 test pairs.

We compare our methods to several conventional gradient-based and neural operator learning baselines as described below: 
\begin{itemize}
    \item $\text{L-BFGS}_{2500}$: Standard L-BFGS optimization to \eqref{eqn:rec} with a maximum of 2500 iterations.
    \item $\text{L-BFGS}_{2500} + \ell_2$ regularization:
    L-BFGS with a maximum of 2500 iterations to \eqref{eqn:rec}, with an additional penalty term $\lambda \|\sigma\|^2$ in the loss function. The regularization coefficient $\lambda$ is data-dependent and was selected from the set $\{10^{-i}\}_{i=3}^9$ that yields the best performance.
    \item ResNet: 
    A ResNet of exactly the same architecture as in Section \ref{sec:ResNet} was trained over the datasets $\{(M_i,\sigma_i^s)\}_{i=1}^N$ to directly learn the target operator $\mathcal{F}: M \mapsto \sigma$.
    \item DDPM: A conditional DDPM model as described in Section \ref{sec:cDDPM} was trained to generate samples of the posterior distribution $p(\sigma \ |\ M)$ over the datasets $\{(M_i,\sigma_i^s)\}_{i=1}^N$.
    A mean estimator over $10$ posterior samples was used as the output. The same estimator is also used for the proposed $\text{L-BFGS}_{K}$ + DDPM model.
\end{itemize}

The proposed methods as well as the baseline models are assessed in terms of the following error metrics.
\begin{itemize}
    \item Relative $\ell_2$ error of the measurement data, which quantifies pixel-wise accuracy on the reconstructed boundary.
    \item Relative $\ell_1$ error of the reconstructed solution, which assesses the pixel-wise accuracy of structural features over the entire image.
    \item Peak-Signal-Noise-Ratio (PSNR), which measures image reconstruction quality by computing the logarithm of the ratio between the maximal value of a signal versus root mean squared error.
    \item Structural Similarity Index Measure (SSIM), which evaluates the similarity of our outputted samples to the distribution of ground truth images by comparing luminance contrast, and structural information.
\end{itemize}

\subsection{Complexity Analysis}

\begin{table}[htb]
    \fontsize{10pt}{12pt}\selectfont
    \centering
    \setlength{\tabcolsep}{4pt}
    \centering
    \begin{tabular}{lrrr}
        \hline
        Model     & Parameters &  Inference Time (s) & Time Complexity \\
        \hline
        $\text{L-BFGS}_{N_{\text{iter}}}$ &  N/A & 172.61 & $O(N_{\text{iter}}n_e n_b)$ \\
        $\text{L-BFGS}_{N_{\text{iter}}}$ + $\ell_2$    & N/A & 167.68 & $O(N_{\text{iter}}n_en_b)$ \\
        ResNet\footnotemark[2]       & 13,271,488 & 0.22 &  $O(n_e)$ \\
        DDPM\footnotemark[2]       & 16,011,265 & 2.11 &  $O(n_e)$   \\
        \midrule
        $\text{L-BFGS}_{K}$ + ResNet  &  13,271,488 & 22.61 & $O(Kn_en_b + n_e)$\\
        $\text{L-BFGS}_{K}$ + DDPM   &  16,011,265 & 27.20 & $O(Kn_en_b + n_e)$ \\
        \hline
    \end{tabular}
    \vspace{0.2 cm}
    \caption{This table presents the number of trainable parameters, the inference time in seconds, and the time complexity during inference with respect to $n_e$ and $n_b$. 
    Our methods are much faster than the conventional L-BFGS baselines, but are slower than vanilla DL methods since we generate initial guesses via $\text{L-BFGS}_{K}$.} 
    \label{tab:time_params_disk}
\end{table}

\footnotetext[2]{For compatibility with the neural network architectures, each D2N measurement $M_i$ is downsampled from its original size $n_b^2$ to $n_e$. Accordingly, the time complexity for ResNet and DDPM is $O(n_e)$ in our implementation.}
In Table~\ref{tab:time_params_disk}, we provide an overview of the number of trainable parameters, inference time and runtime complexity with respect to the original mesh size $n_e$ and number of boundary points $n_b$ for all methods. For all $\text{L-BFGS}_N$ methods, $N$ denotes the maximum number of L-BFGS iterations. L-BFGS methods involve solving the underlying PDE $n_b$ times at each iteration, which has a complexity of $O(n_e n_b)$, while convolutional layers exhibit linear complexity with respect to $n_e$. GPU acceleration enables efficient neural network inference with much faster runtimes compared to standalone L-BFGS methods. We maintain similar parameter counts for ResNet and DDPM to enable a fair comparison in our experiments, although inference time is slower for DDPM due to its sequential denoising paradigm.

Compared to the baselines, our methods achieve $O(Kn_e n_b + n_e)$ complexity, scaling with the number of L-BFGS iterations $K$ as well as the subsequent neural network inference. Compared to standalone L-BFGS methods, our approaches offer substantial speedups of 6-8$\times$ while maintaining enhanced reconstruction quality.

\subsection{Training details}\label{sec:traindets}

A finite element solver is used to solve the elliptic equation \eqref{eqn:elliptic}. For all L-BFGS-utilizing methods, we maintain the number of stored memory updates as $m = 10$.

All deep learning models are trained using the Adam algorithm, with an initial learning rate of $\alpha = 10^{-3}$. 
All DDPM models are trained with $T = 400$ timesteps. We utilize a cosine learning rate scheduler for our DDPM methods as discussed in~\citet{nichol2021cosine}, with a minimum learning rate of $\alpha_{\min} = 10^{-6}$. For training of ResNet models, the learning rate is multiplied by a factor of 0.75 every 500 epochs. No regularization is incorporated in the training of any ResNet models.
All model training took place on an 80 GB NVIDIA A100 GPU. For all datasets, we maintain an image size of $64 \times 64$ pixels and train all deep learning models for 20,000 epochs. We train with batches of size 128 randomly sampled without replacement at each epoch. 

\subsection{Four Circles Dataset}
For each ground truth image $\sigma_i$ from the Four Circles Dataset, we compute the D2N map and solve \eqref{eqn:rec} using L-BFGS for $K = 350$ iterations to obtain a suitable low-fidelity prior $\hat{\sigma}_i$ for the deep learning methods, denoted below as $\text{L-BFGS}_{350}$.

\begin{figure}[htb]
    \centering
    \includegraphics[width=\textwidth]{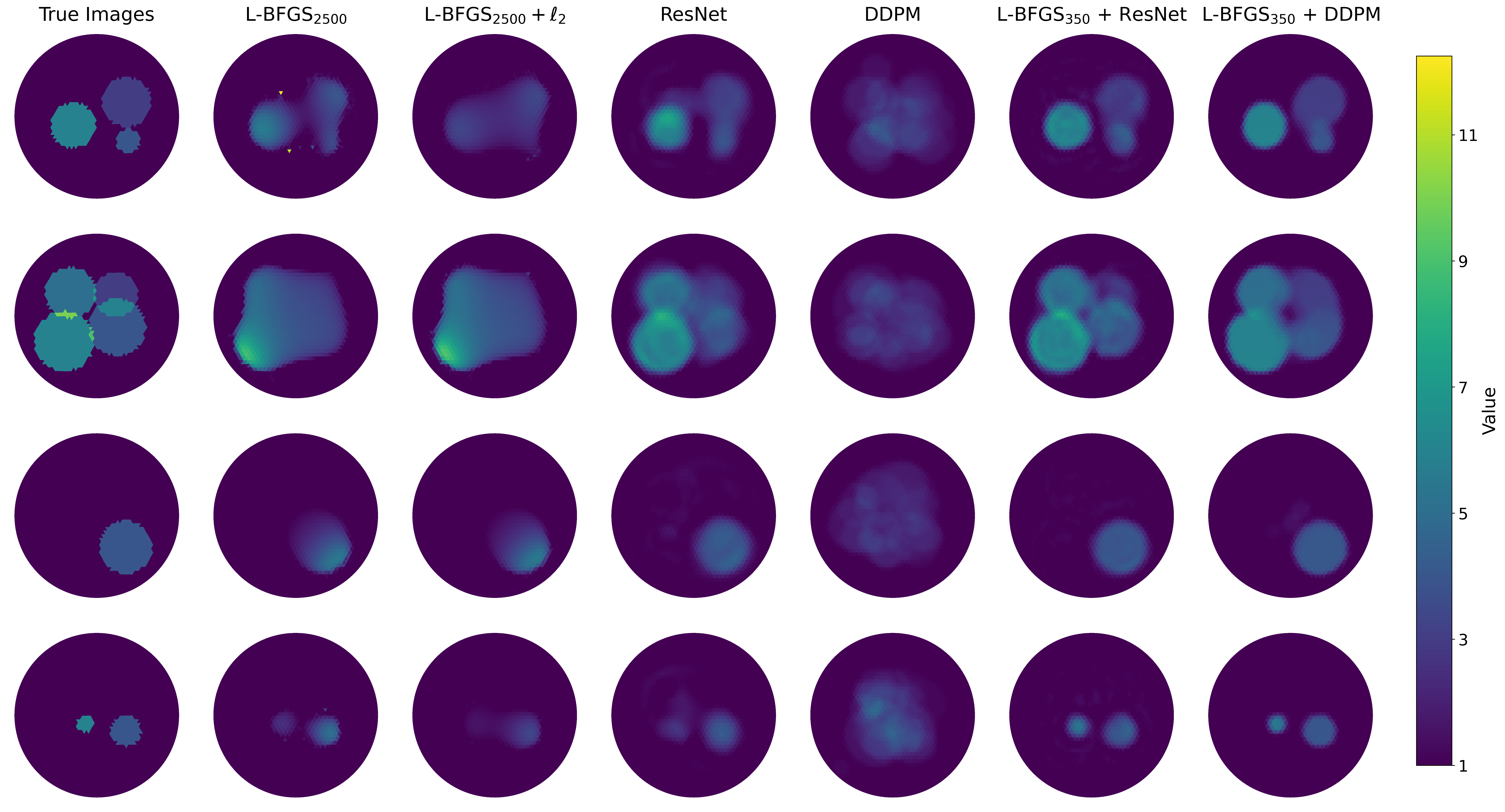}
    \caption{
    \textbf{Four Circles Dataset.} Four different samples of ground truth (Column 1) and reconstructed media
    from baseline models (Columns 2-5) and the proposed methods (Columns 6-7). 
    Our proposed methods perform significantly better than baseline models in capturing the shape and sharp boundary of the circles.
    }
    \label{fig:circs_results_full}
\end{figure}

In Figure~\ref{fig:circs_results_full}, we display the reconstructed images from all four baseline models and our proposed two methods for four different samples.
It is evident that conventional methods perform poorly in determining any structural details of the solution, even with regularization and excessive iterations. 
Additionally, we observe that the ResNet baseline does a decent job in capturing the shape of the circles but introduces blurring artifacts that persist over all numbers of circles. On the other hand, the DDPM baseline, even when averaging to reduce 
sample variance, fails to learn the Four Circles distribution.

In contrast, both our proposed neural correction operator methods demonstrate significant better reconstructions over the baseline models.
In particular, the $\text{L-BFGS}_{350}$ + ResNet method sees sharper boundaries in simpler (1-2 circles) cases. While we still observe persistent blurring for more difficult (3-4 circles) cases, our method is able to capture the shapes of all circles even if some of them are close or intersecting.
Similarly, equipping DDPM with the $\text{L-BFGS}_{350}$ initial guesses results in a dramatic improvement in visual quality from directly learning the target operator.
For 1-3 circles, $\text{L-BFGS}_{350}$ + DDPM obtains the best performance over all methods. While it struggles to capture finer details in the most challenging cases with four intersecting circles, it does a better job generating media with homogeneous regions than ResNet models.

\begin{table}[htb]
\fontsize{9pt}{11.2pt}\selectfont
\centering
\begin{tabular}{l@{\hskip 6pt}c@{\hskip 8pt}c@{\hskip 8pt}c@{\hskip 8pt}c}
\toprule
Model & \parbox[c][0.75cm]{2.4cm}{\centering Rel. $\ell_2$ Error\\ (Measurement)} \hspace{-0.5em} $\downarrow$ & \parbox[c][0.75cm]{1.9cm}{\centering Rel. $\ell_1$ Error\\ (Solution)}\hspace{-0.2em} $\downarrow$ & PSNR $\uparrow$ & SSIM $\uparrow$ \\
\midrule
$\text{L-BFGS}_{2500}$ & $1.0 \times 10^{-4} \pm 2.0 \times 10^{-4}$ & $0.137 \pm 0.064$ & $26.68 \pm 5.67$ & $0.877 \pm 0.069$  \\
${\text{L-BFGS}_{2500} + \ell_2}$ & \textcolor{RoyalBlue}{$1.0 \times 10^{-4} \pm 1.0 \times 10^{-4}$} & $0.155 \pm 0.055$ & $25.36 \pm 3.91$ &	$0.863 \pm 0.060$ \\
ResNet & $0.027 \pm 0.017$ & $0.167 \pm 0.042$ & $27.08 \pm 3.39$ & $0.833 \pm 0.057$  \\
DDPM & $0.021 \pm 0.046$ & $0.429 \pm 0.060$ & $20.06 \pm 2.68$ & $0.654 \pm 0.053$  \\
\midrule
$\text{L-BFGS}_{350}$ + ResNet & $0.029 \pm 0.017$ &  $0.120 \pm 0.040$ & $29.32 \pm 4.07$ & $0.880 \pm 0.053$  \\
$\text{L-BFGS}_{350}$ + DDPM & $0.005 \pm 0.009$ & \textcolor{RoyalBlue}{$0.089 \pm 0.052$} & \textcolor{RoyalBlue}{$29.63 \pm 5.36$} & \textcolor{RoyalBlue}{$0.909 \pm 0.060$}\\
\bottomrule
\end{tabular}
\caption{Mean and standard deviation for error metrics across all models over the \textbf{Four Circles Dataset}. Our methods perform better than all baselines when it comes to metrics that consider accuracy of the sampled images (relative $\ell_1$ error of solution) as well as distribution similarity metrics (PSNR and SSIM).}
\label{tab:circs_stats}
\end{table}

Table~\ref{tab:circs_stats} summarizes the error metrics over the testing data set across all models in the Four Circles distribution. We highlight the mean and standard deviation taken over all images in our test dataset. Our methods perform much better than all baselines when it comes to both pixel-wise accuracy and similarity of reconstructed samples to the original ground truth distribution, which is indicated by higher PSNR and SSIM values. In particular, DDPM benefits the most from using the prior, seeing significant improvement in reducing the relative $\ell_1$ error of the solution, PSNR, and SSIM. While the L-BFGS baselines achieve a lower relative $\ell_2$ error of the measurement data, they fail to capture any visual structure as seen in Figure~\ref{fig:circs_results_full}. This is a manifestation of the ill-posedness of EIT. Incorporating the initial estimates from L-BFGS into deep learning models helps to mitigate this and yields substantial improvements in reconstruction quality while also maintaining measurement accuracy.

\subsection{Shepp-Logan Dataset}

Next, we proceeded to test the neural correction operator using the Shepp-Logan distribution. Here, we run the L-BFGS solver for $K = 150$ iterations (denoted $\text{L-BFGS}_{150}$) to obtain the low-fidelity initial guesses for the deep learning methods. 

\begin{figure}[htb]
    \centering
    \includegraphics[width=0.95\linewidth]{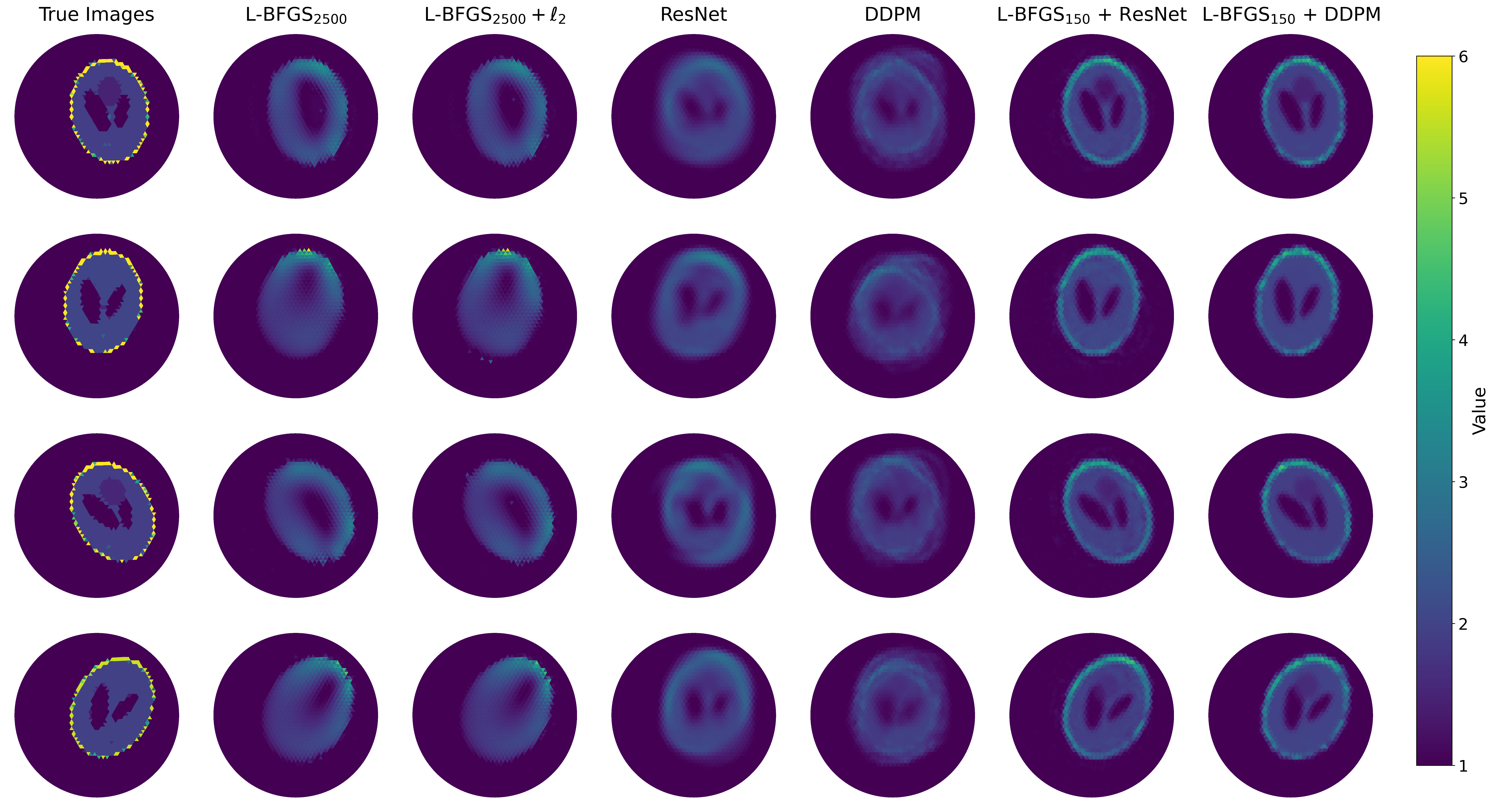}    
    \caption{
    \textbf{Shepp-Logan Dataset.} Four different samples of Ground truth (Column 1) and reconstructed media
    from baseline models (Columns 2-5) and the proposed methods (Columns 6-7). 
    Our proposed methods perform significantly better than baseline models in finding the interior structures.}
    \label{fig:sl_results_full}
\end{figure}

In Figure~\ref{fig:sl_results_full}, we display four samples of ground truth media and the outputs of each method. We observe that the conventional L-BFGS baselines can recover the shape of the largest ellipse but fail to resolve any interior structure. ResNet and DDPM are able to reconstruct the interior slightly, but fail to capture the overall shape and/or orientation of the media. In particular, none of the baseline methods are capable of accurately determining any of the interior ellipses, due to the severe ill-posedness of the EIT problem.

In contrast, neural correction operator methods are able to capture the interior of the media quite well with regards to both shapes and pixel values of the interior ellipses. However, we note some differences in behavior between our two approaches, the $\text{L-BFGS}_{150}$ + ResNet method and the $\text{L-BFGS}_{150}$ + DDPM method. 
We observe that the baseline ResNet method displays a significant amount of blurring within the interior of the media. $\text{L-BFGS}_{150}$ + ResNet mitigates this blurring of the interior, however some slight blurring can still be observed.
On the other hand, $\text{L-BFGS}_{150}$ + DDPM does not admit any blurring artifacts, which is helped by taking the displayed sample as an average over 10 images. 

Additionally, we notice that $\text{L-BFGS}_{150}$ + ResNet has a tendency to produce overly smooth boundaries for the Shepp-Logan phantoms, considering the incomplete boundary of the ground truth Shepp-Logan example in Row 2 of Figure~\ref{fig:sl_results_full}. 
The inability of $\text{L-BFGS}_{150}$ + ResNet to accurately reflect these discontinuities is likely due to the architecture of ResNet, which consists of convolutional layers without upscaling, leading to a preference for capturing a smooth and continuous boundary. In contrast, $\text{L-BFGS}_{150}$ + DDPM more accurately captures these discontinuities in addition to learning the structural details of the interior. Notably, both methods sometimes introduce interior features not present in the ground truth, notably the presence or lack of a circle in the center. We discuss these ``hallucinations'' in further detail in Section~\ref{sec:noise}.

\begin{table}[htb]
    \fontsize{9pt}{11.2pt}\selectfont
    \centering
    \begin{tabular}{l@{\hskip 6pt}c@{\hskip 8pt}c@{\hskip 8pt}c@{\hskip 8pt}c}
    \toprule
    Model & \parbox[c][0.75cm]{2.4cm}{\centering Rel. $\ell_2$ Error\\ (Measurement)} \hspace{-0.5em} $\downarrow$ & \parbox[c][0.75cm]{1.9cm}{\centering Rel. $\ell_1$ Error\\ (Solution)}\hspace{-0.2em} $\downarrow$ & PSNR $\uparrow$ & SSIM $\uparrow$ \\
    \midrule
    $\text{L-BFGS}_{2500}$ & \textcolor{RoyalBlue}{$4.7 \times 10^{-5} \pm 4.3 \times 10^{-5}$} & $0.156 \pm 0.007$ & $19.71 \pm 0.43$ & $0.769 \pm 0.015$  \\
    $\text{L-BFGS}_{2500}$ + $\ell_2$ & $1.0 \times 10^{-4} \pm 1.0 \times 10^{-5}$ & $0.154 \pm 0.008$ & $19.69 \pm 0.42$ & $0.768 \pm 0.014$  \\
    ResNet & $0.009 \pm 0.003$ & $0.193 \pm 0.010$ & $19.08 \pm 0.48$ & $0.696 \pm 0.019$  \\
    DDPM & $0.009 \pm 0.017$ & $0.202 \pm 0.017$ & $18.39 \pm 0.50$ & $0.672 \pm 0.030$  \\
    \midrule
    $\text{L-BFGS}_{150}$ + ResNet & $0.027 \pm 0.006$ & $0.138 \pm 0.006$ & $21.29 \pm 0.42$ & $0.816 \pm 0.010$  \\
    $\text{L-BFGS}_{150}$ + DDPM & $0.005 \pm 0.004$ & \textcolor{RoyalBlue}{$0.123 \pm 0.006$} & \textcolor{RoyalBlue}{$21.39 \pm 0.45$} & \textcolor{RoyalBlue}{$0.829 \pm 0.010$} \\
    \bottomrule
    \end{tabular}
    \caption{Mean and standard deviation for all models over the \textbf{Shepp-Logan Dataset}. Our methods perform better than all baselines when it comes to metrics that consider accuracy of the sampled images (relative $\ell_1$ error of solution) as well as distribution similarity metrics (PSNR and SSIM).}
    \label{tab:sl_stats}
\end{table}

Table~\ref{tab:sl_stats} summarizes the error metrics across all models over the testing dataset for the Shepp-Logan distribution. We note that traditional L-BFGS methods consistently outperform naive deep learning approaches, which do far poorer in simply learning the shape of the distribution. Equipping ResNet and DDPM with the L-BFGS initial guess leads to significant improvements in all metrics, including improvements in the solution $\ell_1$ error of $28\%$ and $39\%$ for ResNet and DDPM respectively. The majority of our error stems from challenges in accurately capturing the values on the boundary, which is a different scale than the rest of the interior (as opposed to the Four Circles distribution). 
Boundary estimation is a common issue when considering the EIT inverse problem~\citep{chen2024pseudo}. Despite these limitations, we are otherwise able to successfully capture all shapes and interior values for a harder distribution than Four Circles.

\subsection{Robustness to Noise}\label{sec:noise}

In this section, we evaluate the robustness of the neural correction operator methods to reconstruct Shepp-Logan images under noisy measurements. To this end, we corrupt the D2N measurements with $1\%$ multiplicative noise $\eta \sim U[-0.01, 0.01]$ and $5\%$ multiplicative noise $\eta \sim U[-0.05, 0.05]$.
For each noisy measurement $M^\eta$, we then run the L-BFGS solver for $K = 250$ iterations (denoted $\text{L-BFGS}_{250}$) to obtain a suitable low-fidelity prior $\hat{\sigma}$ for our methods. 
All models in this subsection were trained and tested with noisy data. 

\begin{figure}[htb]
    \centering
    \includegraphics[width=0.95\textwidth]{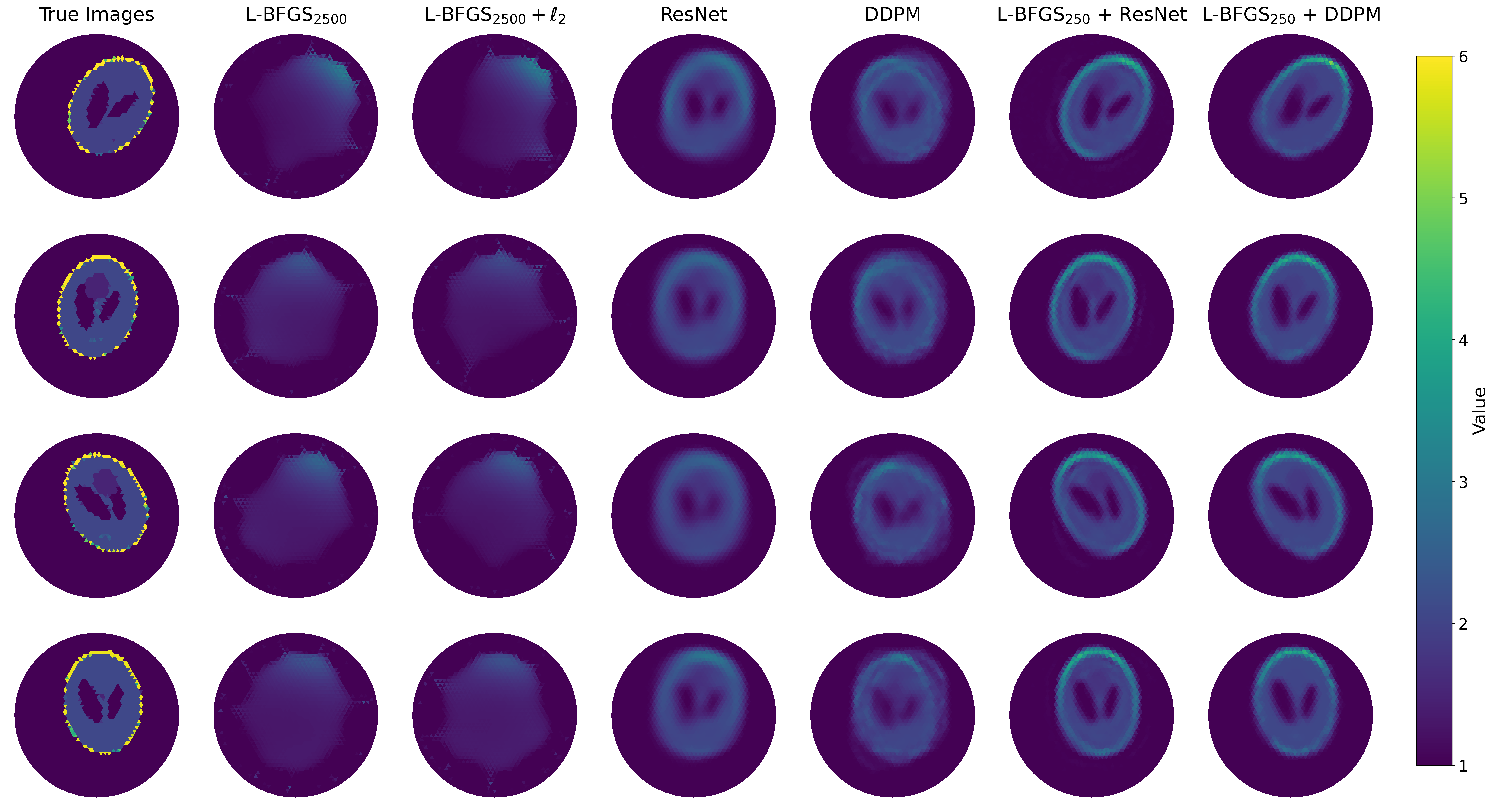}
    \caption{
    \textbf{1\% Noisy Shepp-Logan Dataset.} Four different samples of ground truth (Column 1) and reconstructed media
    from baseline models (Columns 2-5) and the proposed methods (Columns 6-7). 
    Our proposed methods can still capture the overall shape of interior structures. However, $\text{L-BFGS}_{250}$ + DDPM starts to hallucinate with missing interior ellipses.}
    \label{fig:sl_noise_results_full}
\end{figure}

In Figure~\ref{fig:sl_noise_results_full}, we display four samples of Shepp-Logan media and the reconstruction results of each method from measurements corrupted with 1\% noise. The L-BFGS baselines fail to recover any meaningful structure, while both ResNet and DDPM produce samples that are visually comparable to the noiseless case. These models are still able to capture the overall shape of the image, but miss finer details in the interior. 
In comparison, our neural correction operator methods demonstrate robust construction from noisy measurements in terms of both pixel values and shape of the interior ellipses. 
We also measure error statistics over the testing data set in Table \ref{tab:sl_noise_stats}, where it is shown that our methods achieve the best performance comparing to all baselines.

\begin{table}[htb]
\fontsize{9pt}{11.2pt}\selectfont
\centering
\begin{tabular}{l@{\hskip 6pt}c@{\hskip 8pt}c@{\hskip 8pt}c@{\hskip 8pt}c}
\toprule
Model & \parbox[c][0.75cm]{2.4cm}{\centering Rel. $\ell_2$ Error\\ (Measurement)} \hspace{-0.5em} $\downarrow$ & \parbox[c][0.75cm]{1.9cm}{\centering Rel. $\ell_1$ Error\\ (Solution)}\hspace{-0.2em} $\downarrow$ & PSNR $\uparrow$ & SSIM $\uparrow$ \\
\midrule
$\text{L-BFGS}_{2500}$ & $1.0 \times 10^{-2} \pm 9.0 \times 10^{-6}$ & $0.219 \pm 0.008$ & $18.42 \pm 0.37$ & $0.672 \pm 0.010$  \\
$\text{L-BFGS}_{2500}$ + $\ell_2$ & $1.0 \times 10^{-2} \pm 1.1 \times 10^{-5}$ & $0.219 \pm 0.008$ & $18.41 \pm 0.36$ & $0.671 \pm 0.011$  \\
ResNet & $0.008 \pm 0.004$ & $0.198 \pm 0.011$ & $18.92 \pm 0.48$ & $0.688 \pm 0.023$  \\
DDPM & $0.009 \pm 0.014$ & $0.210 \pm 0.016$ & $18.36 \pm 0.49$ & $0.663 \pm 0.030$  \\
\midrule
$\text{L-BFGS}_{250}$ + ResNet & $0.020 \pm 0.011$ & \textcolor{RoyalBlue}{$0.157 \pm 0.008$} & \textcolor{RoyalBlue}{$20.51 \pm 0.44$} & \textcolor{RoyalBlue}{$0.777 \pm 0.013$} \\
$\text{L-BFGS}_{250}$ + DDPM & \textcolor{RoyalBlue}{$0.006 \pm 0.005$} & $0.160 \pm 0.014$ & $19.86 \pm 0.65$ & $0.765 \pm 0.028$ \\
\bottomrule
\end{tabular}
\caption{Error metrics across all models over the Shepp-Logan distribution with 1\% multiplicative noise added to the measurement data. 
In the noisy measurement setting, ResNet as the correction model demonstrates better reconstruction quality compared to DDPM with respect to both pixel-wise accuracy (relative $\ell_1$ error of solution) as well as distribution similarity (PSNR and SSIM).}
\label{tab:sl_noise_stats}
\end{table}

\begin{figure}[htb]
    \centering
    \includegraphics[width=0.95\textwidth]{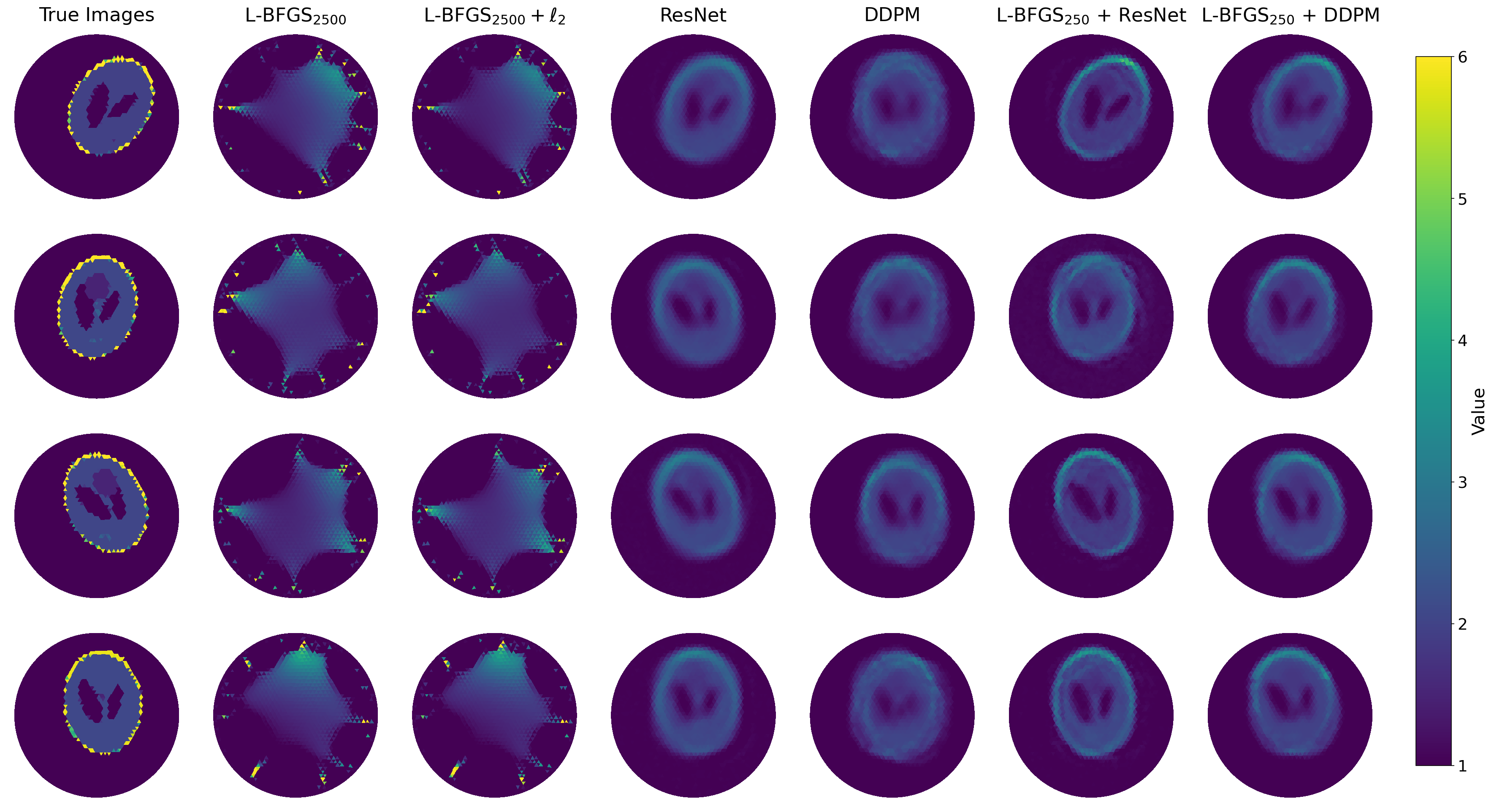}
    \caption{
    \textbf{5\% Noisy Shepp-Logan Dataset.} Four different samples of ground truth (Column 1) and reconstructed media
    from baseline models (Columns 2-5) and the proposed methods (Columns 6-7). 
    At higher noise levels, although our proposed methods still outperform all baselines, they begin to exhibit errors in predicting the boundary and orientation of the media.}
    \label{fig:sl_noise_5pct_results_full}
\end{figure}

Figure~\ref{fig:sl_noise_5pct_results_full} presents the reconstruction results for the same four samples of Shepp-Logan media when the measurements are corrupted with 5\% noise. At this higher noise level, we observe that $\text{L-BFGS}_{250}$ + ResNet displays increased blurring. Additionally, $\text{L-BFGS}_{250}$ + ResNet now struggles to determine the rotational orientation of the media. 
Same error statistics are conducted to the testing dataset with 5\% noise in Table \ref{tab:sl_noise_5pct_stats}, our methods still outperform all baseline models.

While $\text{L-BFGS}_{250}$ + DDPM yields drastic improvements compared to DDPM alone, we observe several limitations of our method under noisy conditions.
$\text{L-BFGS}_{250}$ + DDPM becomes less accurate at reconstructing the overall shape of the image and displays inconsistency when reconstructing the interior features. In particular, the method fails to capture interior artifacts of the media when it is present and/or hallucinates them when they are not there. This is exemplified by the small blue circle in the interior of the media in Rows 2 and 3 of Figures~\ref{fig:sl_noise_results_full} and~\ref{fig:sl_noise_5pct_results_full}. This is likely due to the increased variance in the posterior distribution caused by the noise introduced into the measurements. In the noiseless case as highlighted in Figure~\ref{fig:sl_results_full}, we are able to accurately reconstruct these circles as well as obtain the correct overall shape of the media.

\begin{figure}[htp]
    \centering
    \includegraphics[width=0.95\linewidth]{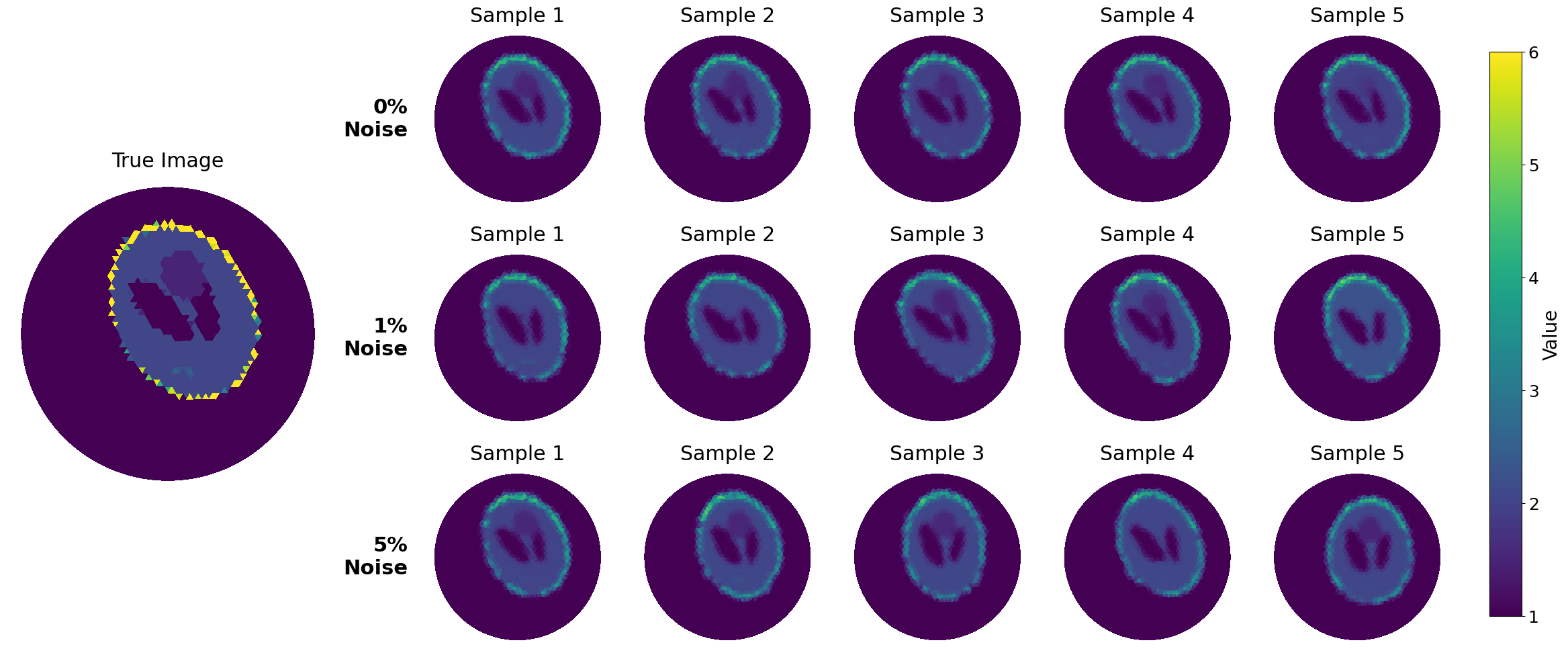}
    \caption{We display five individual $\text{L-BFGS}_{250}$ + DDPM results without averaging for the 0\%, 1\%, and 5\% noise levels. As the noise level increases, the variance of individual samples significantly increases and the accuracy of predicted internal features as well as overall shapes decays.}
    \label{fig:sl_noise_DDPM_comp}
\end{figure}

\begin{table}[htp]
\fontsize{9pt}{11.2pt}\selectfont
\centering
\begin{tabular}{l@{\hskip 6pt}c@{\hskip 8pt}c@{\hskip 8pt}c@{\hskip 8pt}c}
\toprule
Model & \parbox[c][0.75cm]{2.4cm}{\centering Rel. $\ell_2$ Error\\ (Measurement)} \hspace{-0.5em} $\downarrow$ & \parbox[c][0.75cm]{1.9cm}{\centering Rel. $\ell_1$ Error\\ (Solution)}\hspace{-0.2em} $\downarrow$ & PSNR $\uparrow$ & SSIM $\uparrow$ \\
\midrule
$\text{L-BFGS}_{2500}$ & $0.029 \pm 0.002$ & $0.271 \pm 0.021$ & $16.22 \pm 0.48$ & $0.553 \pm 0.036$  \\
$\text{L-BFGS}_{2500}$ + $\ell_2$ & $0.028 \pm 0.002$ & $0.272 \pm 0.021$ & $16.18 \pm 0.48$ & $0.552 \pm 0.036$  \\
ResNet & $0.021 \pm 0.008$ & $0.194 \pm 0.016$ & $19.03 \pm 0.60$ & $0.699 \pm 0.029$  \\
DDPM & $0.006 \pm 0.002$ & $0.207 \pm 0.015$ & $18.48 \pm 0.51$ & $0.662 \pm 0.029$  \\
\midrule
$\text{L-BFGS}_{250}$ + ResNet & $0.042 \pm 0.019$ & \textcolor{RoyalBlue}{$0.179 \pm 0.024$} & \textcolor{RoyalBlue}{$19.72 \pm 0.74$} & \textcolor{RoyalBlue}{$0.740 \pm 0.035$} \\
$\text{L-BFGS}_{250}$ + DDPM & \textcolor{RoyalBlue}{$0.006 \pm 0.003$} & $0.179 \pm 0.018$ & $19.15 \pm 0.68$ & $0.721 \pm 0.035$ \\
\bottomrule
\end{tabular}
\caption{Error metrics across all models over the Shepp-Logan distribution with 5\% multiplicative noise added to the measurement data. 
In the noisy measurement setting, ResNet as the correction model demonstrates better reconstruction quality compared to DDPM with respect to both pixel-wise accuracy (relative $\ell_1$ error of solution) as well as distribution similarity (PSNR and SSIM).}
\label{tab:sl_noise_5pct_stats}
\end{table}

To inspect the degradation in sample quality, we consider the ground truth displayed in Row 3 of Figures~\ref{fig:sl_noise_results_full} and~\ref{fig:sl_noise_5pct_results_full}. 
We display five samples each from $\text{L-BFGS}_{150}$ + DDPM in the noiseless measurement setting and by $\text{L-BFGS}_{250}$ + DDPM with 1\% and 5\% noisy measurement data in Figure~\ref{fig:sl_noise_DDPM_comp}.

In the noiseless setting, $\text{L-BFGS}_{150}$ + DDPM correctly reconstructs the shape of the ground truth in all samples and captures the interior circle in four out of five examples.
However, under 1\% measurement noise, $\text{L-BFGS}_{250}$ + DDPM can still learn the shapes of the two interior ellipses but fails to learn the shape of the overall media in almost all of the displayed samples. In addition to this, we note that $\text{L-BFGS}_{250}$ + DDPM only captures the interior circle in two out of the give displayed samples in Figure~\ref{fig:sl_noise_DDPM_comp}. When the measurement is corrupted with 5\% noise, we observe that in addition to the failure to capture the shape of the media, $\text{L-BFGS}_{250}$ + DDPM fails to learn the orientation of the media, as evidenced by Samples 3 and 5 of Figure~\ref{fig:sl_noise_DDPM_comp}. We believe this is a manifestation of the severe illposedness of EIT.

Table~\ref{tab:noise_comp} displays the relative $\ell_1$ error of all pictured samples and highlights the degradation in performance across the 0\%, 1\% and 5\% noise settings. We observe that while mean reconstruction error increases by $31\%$ between the 0\% and 5\% cases, the standard deviation of the error increases by five times. It is evident that DDPM exhibits far more variance in sample quality when measurement noise is introduced.

\begin{table}[htb]
    \fontsize{8pt}{10pt}\selectfont
    \centering
    \begin{tabular}{c*{7}{S[table-format=1.2, round-mode=places, round-precision=3]}}
    \toprule
    \multirow{3}{*}{Noise Level} & \multicolumn{5}{c}{Relative $\ell_1$ Error of Solution} & \multicolumn{2}{c}{Statistics} \\
    \cmidrule(lr){2-6} \cmidrule(lr){7-8}
    & {Sample 1} & {Sample 2} & {Sample 3} & {Sample 4} & {Sample 5} & {Mean} & {Std. Dev.} \\
    \midrule
    0\% & 0.130858 & 0.122175 & 0.129296 & 0.127779 & 0.135952 & 0.129212 & 0.004077 \\
    1\% & 0.147955 & 0.170394 & 0.154623 & 0.141805 & 0.164438 & 0.155843 & 0.009543 \\
    5\% & 0.1293 & 0.166395 & 0.175029 & 0.168112 & 0.216446 & 0.171056 & 0.025331 \\
    \bottomrule
    \end{tabular}
    \caption{Relative $\ell_1$ error of solution for all samples displayed in Figure~\ref{fig:sl_noise_DDPM_comp}. }
  \label{tab:noise_comp}
\end{table}

\section{Discussion and Conclusion} \label{sec:conclusion}

We proposed the neural correction operator method for solving EIT inverse problems that combines the L-BFGS method with neural network-based operators, including ResNet and conditional diffusion models. 
Our results demonstrate that this decomposition strategy significantly improves reconstruction quality over both standalone optimization methods and direct operator learning approaches. 
Despite the promise of deep learning in scientific computing, we highlight a critical limitation of operator learning in PDE inverse problems: the inherent ill-posedness of EIT leads to a hallucination effect in data-driven models, where plausible-looking but incorrect reconstructions can occur, a challenge underscored in recent work on linear inverse problems~\citep{colbrook2022difficulty}. Addressing these issues will be essential for deploying AI-driven inverse solvers in high-stakes applications such as medical imaging. 
One promising direction for future work is to combine the proposed approach with the warm-initiated methods~\citet{zhou2023neural,guo2025warm} with subsequent gradient-based iterations, which may lead to finer improvements in the reconstructed images.

\section*{Acknowledgments}
 C. W. was partially supported by the National Science Foundation under award DMS-2206332.  
 K. C. was partially supported by the startup fund from the University of Delaware.
\pagebreak

%%       \citet{<label>} ==> Lamport (1994)
%%       \citep{<label>} ==> (Lamport, 1994)
%%
%% Example citation, See \citet{lamport94}.





% %% For authoryear reference style

% \bibitem[Lamport(1994)]{lamport94}
%   \textit{\LaTeX: a document preparation system},
%   Addison Wesley, Massachusetts,

% \end{thebibliography}

\end{document}